\documentclass[dvips,noinfoline,ejs]{imsart}

\doi{10.1214/07-EJS131}
\pubyear{2008}
\volume{2}
\firstpage{118}
\lastpage{126}

\startlocaldefs

\usepackage[latin1]{inputenc}
\usepackage{amsfonts}
\usepackage{amsmath}
\usepackage{qsymbols}
\usepackage{mathbh}

\RequirePackage[colorlinks]{hyperref}

\def \B{{\bf B}}

\def \b{{\sf b}}

\def \be{\begin{eqnarray*}}
\def \ee{\end{eqnarray*}}
\def \ben{\begin{eqnarray}}
\def \een{\end{eqnarray}}
\def \bq{\begin{equation}}
\def \eq{\end{equation}}

\def \build#1#2#3{\mathrel{\mathop{\kern 0pt#1}\limits_{#2}^{#3}}}
\def \cro#1{\llbracket#1\rrbracket}

\def \floor#1{\lfloor#1\rfloor}
\def\proof{\noindent{\bf Proof. }}

\def \eref#1{(\ref{#1})}
\def \sous#1#2{\mathrel{\mathop{\kern 0pt#1}\limits_{#2}}}
\def \sur#1#2{\mathrel{\mathop{\kern 0pt#1}\limits^{#2}}}

\def \dd{\xrightarrow[n]{(d)}}
\def \proba{\xrightarrow[n]{(proba)}}
\def \as{\xrightarrow[n]{(a.s.)}}

\def \dis{\displaystyle}
\def \tend{\longrightarrow}

\def \1{\mathbb{I}}
\def\l{\left}
\def\r{\right}

\def \bar{\overline}
\def \cov{\textrm{cov}}

\def \ind{\mathbh{1}}

\DeclareMathOperator{\mult}{mult}
\newtheorem{lem}{Lemma}

\newtheorem{pro}[lem]{Proposition}
\newtheorem{theo}[lem]{Theorem}

\newtheorem{remi}{Remark\rm}{\rm}
\newtheorem{note}{Note}{}
{\rm}
{\begin{center}\begin{minipage}{16cm}\begin{remi}}%
{\end{remi}\end{minipage}\end{center}}
\endlocaldefs

\begin{document}

\begin{frontmatter}
\title{One more approach to the convergence of the empirical process to the Brownian bridge}
\runtitle{Convergence of the empirical process}

\begin{aug}
\author{\fnms{Jean-Fran\c{c}ois} \snm{Marckert}\ead[label=e1]{marckert@labri.fr}}
\address{CNRS, LaBRI, Universit\'e Bordeaux 1\\
 351 cours de la Lib\'{e}ration\\
33405 Talence cedex, France}

\runauthor{J.-F. Marckert}
\end{aug}

\begin{abstract}A theorem of Donsker asserts that the empirical process converges in distribution to the Brownian bridge. The aim of this paper is to provide a new and simple proof of this fact.
\end{abstract}

\begin{keyword}[class=AMS]
\kwd[Primary ]{62G30}
\kwd{60F17}
\end{keyword}

\begin{keyword}
\kwd{Empirical process}
\kwd{Donsker Theorem}
\kwd{Brownian bridge}
\end{keyword}

\received{\smonth{10} \syear{2007}}

\end{frontmatter}

Let $(U_i)$ be a sequence of i.i.d. random variables uniformly distributed on $[0,1]$, and let $F_n$ be the so-called cumulative empirical function, associated with the $n$ first $U_i$'s:
\[F_n(t):=n^{-1}\#\l\{U_i\leq t, i\in\{1,\dots,n\}\r\},~~~~t\in[0,1].\]
The sequence of processes $(F_n)$ converges pointwise a.s. on [0,1] to $F$ defined by $F(t)=t$; this is a consequence of the strong law of large numbers. The Glivenko-Cantelli theorem asserts that this a.s. convergence stands also for the uniform convergence: a.s. $\sup_{x\in[0,1]}|F_n(x)-F(x)|\sous{\to}{n} 0$. To see this, take $0=x_1<\dots<x_k=1$ and check that by monotonicity of $F_n$ and $F$, $\sup_{x\in[0,1]}|F_n(x)-F(x)|\leq \max_j \max(|F_n(x_{j+1})-F(x_j)|,|F_n(x_{j})-F(x_{j+1})|)\as \max x_{j+1}-x_{j}$, which may be chosen as small as wanted.

In some sense, Donsker's Theorem \cite{D54} provides the second term in this convergence. Consider
\begin{equation}
\b_n(t):=\sqrt{n}\l(F_n(t)-F(t)\r),~~~ t\in[0,1].
\end{equation}
\begin{theo}\label{donsk} (Donsker \cite{D54}) The sequence $(\b_n)$ converges in distribution to the Brownian bridge $\b$ on $D[0,1]$ the space of c\`{a}dl\`{a}g functions on [0,1] equipped with the Skorohod topology.
\end{theo}

\begin{note}
When the variables $U_i$'s are not uniform, the study of the empirical process reduces to the uniform case thanks to a classical ``time change'' involving the inverse of the cumulative function of the $U_i$'s. Some problems of continuity arise due to the atoms of the $U_i$'s but roughly speaking one may say that all the difficulties are present in the case of the uniform distribution. \medskip
\end{note}

We recall that the Brownian bridge is the continuous centered Gaussian process such that $\cov(\b(s),\b(t))=s(1-t)$ when $0\leq s \leq t\leq 1$. It owns the following trajectorial representation~:
\begin{equation}
(\b(t))_{t\in[0,1]}\sur{=}{(d)}(\B_t-t\,\B_1)_{t\in[0,1]},
\end{equation}
where $\B$ is the standard Brownian motion. This may immediately be checked using that $\B$ is a centered Gaussian process such that $\cov(\B_s,\B_t)=\min(s,t)$.

In fact, Donsker proves only in details $\max \b_n\dd \max \b$ justifying the Doob's heuristic \cite{Doob}.  One may find in the literature numerous more or less direct proofs of Theorem \ref{donsk}. See e.g. Billingsley \cite{B68} (and references therein), Kallenberg \cite{KAL}, and also more advanced proofs  and constructions (and stronger results) as that of Koml\'os, Major and Tusn\'ady \cite{KMT}. Some books are devoted to the convergence of empirical measures and processes~: we send the interested reader to Shorack \& Wellner \cite{SW}, van der Vaart \& Wellner \cite{VW}.
As a matter of fact, usual proofs of Theorem \ref{donsk} use often some advanced constructions or are treated in probability books when a lot of materials have been introduced, leading to some intricate and complex proofs, quite difficult to be taught entirely to beginners. The aim of this paper is to present a new proof of Theorem \ref{donsk} using only ``simple'' arguments: only immediate considerations about the weak convergence in $C[0,1]$ and $D[0,1]$, and the other very famous Donsker's Theorem which says that a rescaled random walk converges to the Brownian motion are used. The Appendix recalls this material. \medskip

\indent We begin the proof of Theorem \ref{donsk} following the steps of Donsker \cite{D54}.
We say that a random vector $(M_i)_{i=1,\dots, n}$ has the multinomial distribution with parameters $(k,p_1,\dots,p_n)$, we write $(M_i)_{i=1,\dots, n}\sim \mult(k,p_1,\dots,p_n)$, when
$\mathbb{P}(M_i=m_i,i=1,\dots,n)=\frac{k!}{\prod_{i=1}^n m_i!}\prod_{i=1}^k p_i^{m_i}$
for any prescribed non negative integers $m_1,\dots,m_k$ summing to $k$, and 0 otherwise. \par
The vector $(N_j)_{j=1,\dots,n}$ defined by
\[N_j:=\#\l\{i \in\{1,\dots,n\}, U_i\in \l[(j-1)/n, j/n\r]\r\},\]
 has the $\mult(n,1/n,\dots,1/n)$ distribution.
The empirical process taken at time $k/n$ for $k\in\{0,\dots,n\}$ is a simple function of this vector~:
\begin{equation}\label{kow}
\b_n(k/n)=\sqrt{n}\l(F_n(k/n)-F(k/n)\r)=n^{-1/2}\sum_{j=1}^k (N_j-1).
\end{equation}
Let $\bar{\b_n}$ be the process obtained by interpolating $\b_n$ between the points $\{k/n, k\in\{0,\dots,n\}\}$.

Let $(P_k)$ be a sequence of i.d.d. Poisson random variables with parameter 1. The distribution of $(P_k)_{k=1,\dots,n}$ under the condition $\sum_{k=1}^n P_k=n$ (or $\sum_{k=1}^n (P_k-1)=0$) is also  $\mult(n,1/n,\dots,1/n)$
as can be straightforwardly checked. For any $k\in\{0,\dots,n\}$, set
\begin{equation}\label{poissonrw}
{\bf S}_k=\sum_{j=1}^k (P_j-1)
\end{equation}
and let ${\bf S}=({\bf S}_k)_{k=0,\dots,n}$ be the ``centered'' Poisson random walk, interpolated between integer points. Hence, we have
\begin{equation}\label{ette}
 (\bar{\b_n}(t))_{t\in[0,1]}\sur{=}{(d)}\l(n^{-1/2}{\bf S}_{nt}\r)_{t\in[0,1]}  \textrm{ conditioned by }{\bf S}_n=0,
\end{equation}
and
\begin{equation}\label{et}
\sup_{t\in[0,1]} \l|\bar{\b_n}(t)-\b_n(t)\r|\leq n^{-1/2}\max_{k=1,\dots,n} {N_k}.
\end{equation}
This is controlled as follow: the $N_i$'s are $\textrm{Binomial}(n,1/n)$. By the Markov inequality write
$\mathbb{P}(\max_k N_k\geq \varepsilon\sqrt{n})\leq n\mathbb{P}(N_1\geq \varepsilon\sqrt{n})\leq n \mathbb{E}(e^{N_1})e^{-\varepsilon \sqrt{n}}
=n(1+\frac{e-1}{n})^ne^{-\varepsilon\sqrt{n}}\sim ne^{e-1} e^{-\varepsilon\sqrt{n}}\sous{\to}{n}0$, and then
\begin{equation}\label{rert}
n^{-1/2}\max_{k=1,\dots,n} {N_k}\proba 0.
\end{equation}
Hence, by \eref{ette}, \eref{et} and \eref{rert}, (see also Lemma \ref{zou} in Appendix) Theorem \ref{donsk} stating the convergence of $(\b_n)$ to $\b$ in $D[0,1]$ is easily implied by the following proposition.
\begin{pro}\label{nwp}
The sequence $\l(n^{-1/2}{\bf S}_{nt}\r)_{t\in[0,1]}$ conditioned by ${\bf S}_n=0$
converges in distribution to $\b$ in $C[0,1]$ equipped with the topology of the uniform convergence.
\end{pro}
The proof we propose for this classical proposition is the real novelty of this paper.

\subsection*{The ``correction'' of a Poisson random walk}
 The main line in our approach is the comparison between ${\bf S}$ and ${\bf S}$ conditioned by ${\bf S}_n=0$. We introduce a correcting process ${\bf C}=({\bf C}_k)_{k=0,\dots,n}$, such that the pair  $({\bf S},{\bf C})$ have the following feature~:\\
$\bullet$  ${\bf S}$ is the centered Poisson random walk (defined in \eref{poissonrw}), \\
$\bullet$ ${\bf S}+{\bf C}$ is distributed as ${\bf S}$ conditioned by ${\bf S}_n=0$. \par

\begin{note} Transforming a problem involving $n$ random variables into a problem involving $I_n\sim$ Poisson($n$) random variables is called Poissonization. Taking $U_1,\dots, U_{I_n}$ instead of $U_1,\dots,U_n$ in the construction presented at the beginning of the paper amounts to replacing ${\bf S}$ conditioned by ${\bf S}_n=0$ by the centered Poisson random walk ${\bf S}$. This is a Poissonization. The correction of the Poisson random walk we propose, which allows to pass from ${\bf S}$ to ${\bf S}$ conditioned by ${\bf S}_n=0$ is from our point of view different in nature from the usual depoissonization techniques. Here, everything relies on an exact combinatoral correction, when usually, most rely on the convergence in distribution of $(I_n-n)/\sqrt{n}$, ensuring the problem with $n$ variables and with $I_n$ variables being asymptotically equivalent, which is not the case here.
\end{note}

Let us come back to our correction procedure. To fix the details, we will use a classical interpretation of Poisson random walk in term of urns/balls. \\
Conditionally on ${\bf S}_n=s$, the vector $(P_i)_{i=1,\dots,n}$ has the $\mult(s+n,1/n,\dots,1/n)$ law. When $m$  balls labeled $1,\dots,m$ are sent independently in $n$ urns according to the uniform distribution, the vector $(N'_i)_{i=1,\dots,n}$ giving the number of balls in the urns follows also the $\mult(m,1/n,\dots,1/n)$ distribution. \par
Let us throw  $P_i$ balls in urn $i$ where $(P_i)_{i=1,\dots,n}$ are i.i.d. Poisson random variables with parameter 1. Then three cases arise: $\sum_{i=1}^n P_i=n$ or $\sum_{i=1}^n P_i<n$, or $\sum_{i=1}^n P_i>n$ (recall that ${\bf S}_n=\sum_{i=1}^n P_i-n$).\\
In the first case ${\bf S}_n=0$, no correction are necessary, then set ${\bf C}_i=0$ for any $i$. The two last cases are treated below.
Notice that we focus on the uni-dimensional distributions of the process ${\bf C}$ since this will appear to be sufficient.
\medskip

\noindent{\bf Case ${\bf S}_n<0$. }
We work conditionally on ${\bf S}_n=s$.
Since $-s$ balls are lacking: throw $-s$ new balls and denote by ${\bf C}_k$ the number of new balls fallen in the $k$ first urns; for any $k$,
\begin{equation}\label{eza}
{\bf C}_k \sim \textrm{Binomial}(-s,k/n).
\end{equation}
More precisely, $(\Delta {\bf C}_k)_{k=1,\dots,n}\sim \mult(-s,1/n,\dots,1/n)$ where $\Delta{\bf C}_k:={\bf C}_k-{\bf C}_{k-1}$ is the $k$th increment of the correcting process ${\bf C}$ (with ${\bf C}_0=0$).
\begin{lem}\label{tp}For any $s<0$ and any $n\geq 1$, conditionally on ${\bf S}_n=s$ the process ${\bf S}+{\bf C}$ is distributed as ${\bf S}$ conditioned by ${\bf S}_n=0$ and ${\bf C}_k \sim \textrm{Binomial}(-s,k/n)$.
\end{lem}
\proof  \rm If $X\sim \mult(n+s,1/n,\dots,1/n), Y\sim(-s,1/n,\dots,1/n)$ and $X$ and $Y$ are independent then $X+Y\sim \mult(n,1/n,\dots,1/n)$. ~$\Box$\medskip

\noindent{\bf Case ${\bf S}_n>0$. } We work conditionally on ${\bf S}_n=s$.
In this case $n+s$ balls have been thrown instead of $n$ and then $s$ balls must be taken out. The vector $(V_k)_{k=1,\dots,n}$ giving the exceeding number of balls in the different urns (those with labels in $n+1,\dots,n+s$) follows the law $\mult(s,1/n,\cdots,1/n)$. Then given ${\bf S}_n=s$, we search a correcting process  $({\Delta {\bf C_k}})_{k=1,\dots,n}\sur{=}{(d)}(-V_k)_{k=1,\dots,n}$. Of course there is a problem to define the correcting process in terms of balls/urns, the balls/urns problem living a priori on a larger probability space than the $(P_i)$'s.
But this gives us the intuition for a right correcting process: we define ${\bf C}$  conditionally on the $P_i$'s as follows. Let $(p_i)_{i=1,\dots,n}$ be non negative integers summing to $n+s$. Set
\begin{equation}\label{corre}
\mathbb{P}\l( \Delta{\bf C}_k=-c_k, k\in\{1,\dots,n\} | P_i=p_i, i\in\{1,\dots,n\} \r)=\frac{\prod_{i=1}^n \binom{p_i}{c_i}\ind_{c_i\leq p_i}}{\binom{\sum p_i}{s}}
\end{equation}
for any given  non negative integers $c_1,\dots,c_n$ summing to $s$, and 0 otherwise.

 \begin{lem}\label{tg}For any $s>0$ and any $n\geq 1$, conditionally on ${\bf S}_n=s$ the process ${\bf S}+{\bf C}$ is distributed as ${\bf S}$ conditioned by ${\bf S}_n=0$ and ${\bf C}_k\sim -\textrm{Binomial}(s,k/n)$.
\end{lem}
\sl Proof. \rm We have to check that ${\bf C}+{\bf S}$ is distributed as ${\bf S}$
conditioned by ${\bf S}_n=0$:\\
$\mathbb{P}(P_i+\Delta {\bf C}_i=j_i, \forall i\,|\,{\bf S}_n=s)$
\be
&=& \sum_{(p_i)~:~ \sum p_i=n+s,p_i\geq j_i}\frac{\mathbb{P}\l( \Delta {\bf C}_i=-(p_i-j_i),\forall i | P_i=p_i,\forall i \r)~
\mathbb{P}(P_i=p_i,\forall i)}{\mathbb{P}({\bf S}_n=s)}\\
&=&\sum_{(p_i-j_i)~:~p_i-j_i\geq 0, \sum p_i-j_i=s}
\frac{\prod_{i=1}^n \binom{p_i}{p_i-j_i}}{\binom{n+s}{s}}
\frac{e^{-n}\prod_{i=1}^n \frac1{p_i!}}
{\frac{e^{-n}n^{n+s}}{(n+s)!}}\\
&=&\frac{e^{-n}\prod_{i=1}^n \frac{1}{j_i!}}{e^{-n}n^n/n!}=\mathbb{P}(P_i=j_i, \forall i\,|\,{\bf S}_n=0)
\ee
where we have used \eref{corre}, the fact that $n+{\bf S}_n$ is Poisson$(n)$ distributed, and
\begin{equation}
\label{mutlsum}
\sum_{(\alpha_i)~:~\alpha_i\geq 0, \sum \alpha_i=s}\prod_{i=1}^n\frac{s!}{\alpha_i!}=(1+\dots+1)^s=n^s.
\end{equation}
We now show
that knowing ${\bf S}_n=s$, $(-\Delta {\bf C}_k)_{k=1,\dots,n}\sim \mult(s,1/n,\dots,1/n)$. This implies the second point. Let $c_1,\dots,c_n$ be non negative integers summing to $s$. Write  $\mathbb{P}(\Delta {\bf C}_k=-c_k,\forall k |{\bf S}_n=s)$
\be
&=&\sum_{(p_i), p_i\geq c_i, \sum p_i=n+s} \frac{\mathbb{P}(\Delta {\bf C}_k=-c_k,\forall k |P_i=p_i,\forall i)\mathbb{P}(P_i=p_i,\forall i)}{\mathbb{P}({\bf S}_n=s)}.
\ee
By \eref{mutlsum}, ${\bf S}_n+n\sim$Poisson$(n)$ and \eref{corre}, this is easily shown to be equal to $\frac{s!}{\prod_{i=1}^n {c_i!}}\frac{1}{n^s}$.
~$\Box$\medskip

\noindent From now on, consider the process ${\bf C}$ as being interpolated between integer points.
\begin{lem} \label{compar}
For any $t\in[0,1]$,
\begin{equation}\label{rezze}
n^{-1/2} |{\bf C}_{nt}+t{\bf S}_n|\proba 0.
\end{equation}
\end{lem}
\proof We work with $C_{\floor{nt}}$ instead of $C_{nt}$.
For $t=0$, \eref{rezze} holds clearly. Let $t\in (0,1]$, and $\varepsilon>0$ be fixed.
Write
\[\mathbb{P}(|{\bf C}_{\floor{nt}}+t{\bf S}_n|\geq \varepsilon \sqrt{n})\leq V_n^M+W_n^M\]
 where
\be
V_n^M&=& \mathbb{P}( |{\bf C}_{\floor{nt}}+t{\bf S}_n|\geq \varepsilon\sqrt{n},|{\bf S}_n|\in\sqrt{n}[M^{-1},M]),\\
W_n^M&=&\mathbb{P}(|{\bf S}_n|\notin\sqrt{n}[M^{-1},M]).
\ee
 Let $\alpha>0$ be a fixed (small) positive real number.
The central limit theorem applied to ${\bf S}_n$ ensures that there exists $M$ such that
$W_n^M \leq \alpha,$ for $n$ large enough.  Fix this $M$.
To bound $V_n^M$, use $\mathbb{P}(A| \cup_i B_i)=\sum \mathbb{P}(A |B_i)\mathbb{P}(B_i)/P( \cup B_i)\leq \max_i \mathbb{P}(A|B_i)$ for disjoint sets $B_i$. Hence
\be
V_n^M
&\leq& \max_{k, |k|\in  \sqrt{n}[M^{-1},M]}\mathbb{P}(|{\bf C}_{\floor{nt}}+t{\bf S}_n|\geq \varepsilon\sqrt{n}~|~ |{\bf S}_{n}|=k).\ee
By Lemmas \ref{tp} and \ref{tg}, $\mathbb{P}(|{\bf C}_{\floor{nt}}+t{\bf S}_n|\geq \varepsilon\sqrt{n}~|~ |{\bf S}_{n}|=k)=\mathbb{P}(|B(k,\floor{nt}/n)-tk|\geq \varepsilon\sqrt{n})$ where $B(k,\floor{nt}/n)$ is a binomial random variable with parameters $k$ and  $\floor{nt}/n$.  Further, by the Bienaym\'{e}-Tchebichev inequality
\[ \max_{k, |k|\in  \sqrt{n}[M^{-1},M]} \mathbb{P}(|B(k,\floor{nt}/n)-tk|\geq \varepsilon\sqrt{n})\to 0.~\Box\]

\begin{pro}$(i)$ The following convergence holds in $C([0,1],\mathbb{R}^2)$:
\[n^{-1/2}({\bf S}_{nt},{\bf C}_{nt})_{t\in[0,1]}\dd (\B_t,-t\B_1)_{t\in[0,1]}.\]
$(ii)$  The following convergence holds in $C[0,1]$:
\[n^{-1/2}({\bf S}_{nt}+{\bf C}_{nt})_{t\in[0,1]}\dd (\B_t-t\B_1)_{t\in[0,1]}.\]
\end{pro}
Proposition \ref{nwp} is a consequence of $(ii)$ thanks to Lemmas \ref{tp} and \ref{tg}.

\proof Assertion  $(ii)$ is a consequence of $(i)$. Proof of $(i)$~: the convergence of $n^{-1/2}{\bf S}_{n.}$ to $\B$ in $C[0,1]$ is given by the other famous Donsker's theorem stating the convergence of rescaled random walks to the Brownian motion (see \cite{B68} or \cite{KAL}). In particular
\begin{equation}\label{terpos}
n^{-1/2}{\bf S}_{n}\to \B_1.
\end{equation}
The finite dimensional distribution of $n^{-1/2}{\bf C}_{n.}$  converges to those of the process $(t\B_1)_{t\in[0,1]}$. Indeed by Lemma \ref{compar} and \eref{terpos},
\[n^{-1/2}({\bf S}_n,{\bf C}_{nt_1},\dots,{\bf C}_{nt_k})\dd (\B_1,-t_1 \B_1,\dots,-t_k \B_1)\]
 for any $0\leq t_1\leq \dots \leq t_k\leq 1$. Then the family $(n^{-1/2}{\bf C}_{n.})$ is tight since it is a sequence of monotone processes whose finite dimensional distribution converge to those of the a.s. continuous process $(t\B_1)_{t\in[0,1]}$ (this is Lemma \ref{tret}$(ii)$). Hence the family $n^{-1/2}({\bf S}_{n\cdot},{\bf C}_{n\cdot})$ is tight. The limit is identified again thanks to Lemma~\ref{compar}. ~$\Box$

\subsection*{Conclusion}
\label{aze} The idea of this proof appeared after a discussion with
Philippe Duchon, where he explained me his algorithm to generate
uniformly a Bernoulli bridge with $2n$ steps, that is a random walk
${\bf S}=({\bf S}_k)_{k=0,\dots,2n}$ with increments $\pm1$,
conditioned by ${\bf S}_{2n}=0$~: build first a simple random walk with
$2n$ steps, choosing i.i.d. increments +1,or $-1$ with probability
$1/2$. If ${\bf S}_{2n}=0$ then it's done. If not, assume that ${\bf
S}_{2n}=2k>0$. Then pick up at random indices $I_{1},I_{2}, \dots$ in
$\cro{1,2n}$.  If $I_i$ is the index of a positive increment, change it
into a negative one; if it is negative then do nothing. Stop when you
have changed $k$ increments. By a simple symmetry argument the path
obtained is uniform in the set of Bernoulli bridges of size $2n$. I
found that this was a nice way to prove that rescaled Bernoulli bridge
converges to the Brownian bridge; this can be proved using the same
argument than the ones exposed above: the correction procedure will
asymptotically and ``eventually removes a straight line of the Brownian
motion''. Therefore, I tried to find other increment distributions for
which a similar correction procedure would have been possible.  It
appears to be not so general, or at least, not so agreeable. The
problem is the following one: in general there does not exist any
simple correction procedure that conserves at each step of the
correction the property of the trajectory to have conditionally on its
terminal position $k$, the law of a simple random walk conditioned by
${\bf S}_n=k$.

\section*{Appendix}

\subsection*{A simple link between the convergence in $C[0,1]$ and in $D[0,1]$}

\begin{lem}\label{zou}
 Let $(X_n)$ be a sequence of processes taking their values in $D[0,1]$. Assume that for any $n$,
$X_n=Y_n+Z_n$
where $Y_n$ is a continuous process, and $Z_n$ is a c\`{a}dlag process.
If $(Y_n)$ converges in distribution to $Y$ in $C[0,1]$, and if $\sup|Z_n|\dd 0$ then $(X_n)$ converges in distribution to $Y$ in $D[0,1]$.
\end{lem}
\proof As a matter of fact, this statement is easy if one knows:\\
$(a)$ $(Y_n)\dd Y$ in $C[0,1]$ implies $(Y_n)\dd Y$ in $D[0,1]$,\\
$(b)$ $(Z_n)$ c\`{a}dlag, $\sup|Z_n|\dd 0$ implies $Z_n\dd 0$ (the null process) on $D[0,1]$. \\
Indeed, knowing this, Lemma \ref{zou} is a consequence of the following classical result: let $A_n,B_n,C_n$ be random variables in a metric space $(S,\rho)$. If $A_n=B_n+C_n$, $A_n\dd A$, $\rho(B_n)\dd0$ then $A_n+B_n\dd A$ (see e.g. Billinglsey \cite[Theorem 3.1]{B68}).\par
In order to understand why $(a)$ and $(b)$ hold true, recall that the Skhorohod topology on $D[0,1]$ is defined by a metric: let $\Lambda$ be the set of strictly increasing continuous functions $\lambda$, satisfying $\lambda(0)=0, \lambda(1)=1$. The metric is
\[d_s(f,g)=\inf_{\lambda \in \Lambda}\Big\{\sup\{|\lambda(t)-t|,t\in[0,1]\}  \vee \sup\{|f(\lambda(t))-g(t)|,{t\in[0,1]}\}\Big\}. \]
Hence $d_s(f,g)\leq \sup\{|f(t)-g(t)|,t\in[0,1]\}$; this yields immediately to $(a)$ and~$(b)$.~$\Box$

\subsection*{Convergence in $C[0,1]$ and in $C[0,1]^2$}
We recall some classical facts concerning the weak convergence in $C[0,1]$ and $(C[0,1])^2$. First tightness and relative compactness are equivalent in these sets by Prohorov's theorem, since they are both Polish spaces.
\begin{lem}\label{tret}
$(i)$Let $(X_n,Y_n)$ be a sequence of pairs of processes in  $(C[0,1])^2$. The tightnesses  of both families $(X_n)$ and $(Y_n)$ imply that of $(X_n,Y_n)$.\\
$(ii)$ Let $(X_n)$ be a sequence of monotone processes in $C[0,1]$. If the finite dimensional distributions of $(X_n)$ converge to those of an a.s. continuous process $X$, then $(X_n)$ is tight and then $X_n\dd X$ in $C[0,1]$.
\end{lem}
\proof $(i)$ Take two compacts $K_1$ and $K_2$ of $C[0,1]$ such that $\mathbb{P}(X_n\in K_1)\geq 1-\varepsilon$ and $\mathbb{P}(Y_n\in K_2)\geq 1-\varepsilon$, then $\mathbb{P}((X_n,Y_n)\in K_1\times K_2)\geq 1-2\varepsilon$ and $K_1\times K_2$ is compact in $(C[0,1])^2$.\\
$(ii)$ Only the tightness of $(X_n)$ in $C[0,1]$ has to be checked. For any function $f:[0,1]\to \mathbb{R}$, and $\delta>0$, the global modulus of continuity of $f$ is
\[\omega_\delta(f)=\sup\{|f(x)-f(y)|,x,y\in [0,1],|x-y|\leq \delta\}.\]
Since $X_n$ is increasing, for any positive integer $m$,
\[\omega_{1/m}(X_n)\leq A_{m,n}:=2\max\l\{\l|X_n(\frac{k}m)-X_n(\frac{k-1}{m})\r|, k=1,\dots,m\r\}.\]
Since the finite dimensional distributions of $(X_n)$ converge to those of $X$,
\[A_{m,n}\build{\tend}{n}{(d)} A_m:=2\max\l\{\l|X(\frac{k}m)-X(\frac{k-1}{m})\r|, k=1,\dots,m\r\}\]
and by the uniform continuity of $X$, $A_m\build{\tend}{m}{proba.} 0$.
Hence $\dis\lim_{m} \limsup_n \mathbb{P}(\omega_{1/m}(X_n)\geq
\varepsilon)=0$ for any $\varepsilon>0$.~$\Box$

\subsection*{Acknowledgments} I thanks Philippe Duchon who gives me the inspiration of this approach. I thank also Bernard Bercu for his comments on a preliminary version of this text.

\end{document}